\documentclass[final]{siamltex}

\usepackage{amssymb}
\usepackage{amsmath}
\usepackage{amssymb}
\usepackage{amsbsy}
\usepackage{float}
\usepackage{graphicx}
\usepackage{subfigure}

\newcommand{\norm}[1]{\mathop{\|} #1 \mathop{\|}_2}

\newcommand{\chol}[1]{\mathop{\mathrm{chol}}\left\{#1\right\}}
\newcommand{\V}[1]{\ensuremath{\boldsymbol{#1}}}
\newcommand{\M}[1]{\ensuremath{\mathsf{#1}}}
\newcommand{\VT}[1]{\ensuremath{\boldsymbol{#1}}^{\textrm{T}}}
\newcommand{\MT}[1]{\ensuremath{\mathsf{#1}}^{\textrm{T}}}
\newcommand{\MTI}[1]{\ensuremath{\mathsf{#1}}^{-\textrm{T}}}
\newcommand{\MTb}[1]{\ensuremath{\left(\mathsf{#1}}\right)^{\textrm{T}}}

\newcommand{\me}{\mathrm{e}}

\newcommand{\md}{\,\mathrm{d}}
\newcommand{\defas}{\triangleq}
\usepackage{courier}
\usepackage{xcolor}
\usepackage{listings}
\title{A Matrix Framework for the Solution of ODE\lowercase{s}:\\
Initial-, Boundary-, and Inner-Value Problems}

\author{Matthew Harker\footnotemark[2]
        \and Paul O'Leary\footnotemark[2] }
\begin{document}

\maketitle

\lstset{language=Matlab,basicstyle=\small\ttfamily,
commentstyle=\color[rgb]{.133,.545,.133}\small\ttfamily,
keywordstyle=\color[rgb]{0,0,1}\ttfamily\small,
stringstyle=\color[rgb]{.627,.126,.941}\ttfamily\small,
showstringspaces=false,
xleftmargin=0.8cm,resetmargins=true,escapechar=\$,captionpos=b,
breaklines=true,breakatwhitespace=false,breakindent=2cm
,numbers=left,numbersep=2ex,numberstyle=\small\ttfamily,
columns=flexible}

\renewcommand{\thefootnote}{\fnsymbol{footnote}}

\footnotetext[2]{Institute for Automation, University of Leoben,
Peter Tunner Strasse 27, A8700 Leoben, Austria}
\renewcommand{\thefootnote}{\arabic{footnote}}

\begin{abstract}
A matrix framework is presented for the solution of ODEs,
including initial-, boundary and inner-value problems. The
framework enables the solution of the ODEs for arbitrary nodes.
There are four key issues involved in the formulation of the
framework: the use of a Lanczos process with complete
reorthogonalization for the synthesis of discrete orthonormal
polynomials (DOP) orthogonal over arbitrary nodes within the unit
circle on the complex plane; a consistent definition of a local
differentiating matrix which implements a uniform degree of
approximation over the complete support --- this is particularly
important for initial and boundary value problems; a method of
computing a set of constraints as a constraining matrix and a
method to generate orthonormal admissible functions from the
constraints and a DOP matrix; the formulation of the solution to
the ODEs as a least squares problem. The computation of the
solution is a direct matrix method. The worst case maximum number
of computations required to obtain the solution is known a-priori.
This makes the method, by definition, suitable for real-time
applications.

The functionality of the framework is demonstrated using a
selection of initial value problems, Sturm-Liouville problems and
a classical Engineering boundary value problem. The framework is,
however, generally formulated and is applicable to countless
differential equation problems.
\end{abstract}



\pagestyle{myheadings}
\thispagestyle{plain}
%

\section{Introduction}\label{intro}

%
There are  a number of papers in which the Taylor Matrix is used
to compute solutions to differential
equations~\cite{Kesan2003,Kurt2008}. These methods use the known
analytical relationship between the coefficients $\V{s}$ of a
Taylor polynomial and those of its derivatives $\dot{\V{s}}$ to
compute a differentiating matrix $\M{D}$. The matrix $\M{D}$
together with the matrix of basis functions arranged as the
columns of the matrix $\M{B}$ are used to compute numerical
solutions to the differential equations. The method of the Taylor
matrix was also extended to the computation of fractional
derivatives~\cite{Keskyn2011}. The problem associated with this
approach is that the computation of the numerical solutions
requires the inversion of the Vandermonde matrix, a process which
is known to be numerically unstable, and dependent on the degree
and node placement. The advantage of the Taylor approach lies in
its ability to yield a solution for arbitrary nodes.

A Chebyshev matrix approach was presented by
Sezer~\cite{Sezer1996}. The approach is fundamentally the same as
for the Taylor matrix, whereby the Chebyshev polynomials are used
as an alternative to the geometric polynomials. The main
restriction associated with the Chebvshev polynomial approach is
that the numerical solution to the differential equations is
restricted to the locations of the Chebyshev points; this lacks
the generality needed for many differential equations and
applications.

Podlubny introduced a matrix approach to discrete fractional
calculus~\cite{Podlubny2000} and later extended this work to
partial fractional differential
calculus~\cite{Podlubny2009,Podlubny2009B}. Triangular strip
matrices play a central role in the work; they are used to perform
integration. They implement the integration from a lower to an
upper bound (or vice versa), whereby the errors accumulate as the
integration proceeds. This poses a problem if inverse problems are
addressed, since it gives the solution an implicit direction and a
different accumulation of errors if the problem is solved from
lower to upper bound or from upper to lower. Furthermore, it is
assumed that the initial value is zero. This makes the method
unsuitable for arbitrary boundary conditions. An early source of
this formulation was proposed by Courant et
al.~\cite{Courant1928}, (a later English translation of the paper
is available~\cite{Courant1967}).

A matrix solution specific to Sturm-Liouville problems was
presented by Amodio~\cite{Amodio2011}. The method is specifically
restricted to Sturm-Liouville problems; furthermore, it only
supports solutions on regularly spaced nodes. The results are
correspondingly modest for problems where the Chebyshev points
yield better solutions, e.g., in the solution of the truncated
hydrogen equation. A number of matrix approaches based on the
Numerov method, and modifications of this method, have also been
presented~\cite{Ledoux2007} for the solution of Sturm-Liouville
problems, however, these methods can not be extended to ODEs in
general.

In this paper we formulate a general matrix framework for the
solution of ordinary differential equations, with arbitrary
initial-, boundary-, or inner values. the main contributions of
the paper are:
\begin{enumerate}
    \item The proposal of a consistent framework of matrices and
    solution approaches which can be applied to initial-,
    boundary-, and inner-value problems;
    \item The implementation of new approached to the synthesis of
    discrete orthonormal basis functions, with and without
    weighting;
    \item Generating differentiating
    matrices which are of constant degree of approximation over
    the complete support. It is particularly important that the
    degree of approximation is consistent at the ends of the
    support if initial and boundary value problems are to
    be solved satisfactorily;
    \item The derivation of a means of synthesizing
    constrained basis functions which form orthonormal matrices.
    This basis functions span the space of all solutions which
    fulfil the constraints. They can be used as admissible functions
    in a discrete equivalent of a Rayleigh-Ritz method;
    \item The formulation of the solution of the ODEs as least
    squares approximations. In this manner there is no
    accumulation of errors.
\end{enumerate}

This paper is structured as follows: In Section~(\ref{sec:A}) the
framework for the generation of all the matrices required to
formulate differential equations as matrix linear differential
operators is presented. Section~(\ref{sec:dsolve}) presents the
approach to discretization of the differential equations and their
solution as a least squares minimization is presented. The
required conditions for a unique solution are derived and two
solution approaches are presented: a direct solution in the case
of a unique solution and the implementation of a discrete
Rayleigh-Ritz method for eigenvalue/eigenvector solutions, e.g.,
as encountered in the solution of Sturm-Liouville problems..
Finally, in Section~(\ref{sec:test}) the performance of the
proposed framework is tested with a series of initial-value
problems, Sturm-Liouville problems and a classical Engineering
boundary value problem.
%
\section{Algebraic Framework}\label{sec:A}
In this section we derive the structure and methods for the
synthesis of all matrices required for the discretization and
solution of ordinary differential equations.

\subsection{Quality Measure for Basis Functions}
An objective measure for the quality of a set of basis functions
is required if the sources of numerical error are to be determined
and the best synthesis method is to be selected. In this paper
continuous polynomials are considered which form orthogonal bases
when evaluated over a discrete measure. The basis functions
$\V{b}_i$, i.e., the polynomials evaluated at discrete points, can
be concatenated to form a matrix, $\M{B} = [ \V{b}_1 \ldots
\V{b}_n ]$. The discrete orthogonal polynomials (DOP) are
characterized by the relationship,
\begin{equation}
    \MT{B} \, \M{W} \, \M{B} = \M{I},
\end{equation}
where $\M{W}$ is the weighting matrix. The Gram matrix is defined
as $\M{G}
\defas \MT{B} \, \M{W} \, \M{B}$. Consequently, the orthogonal complement
$\M{G}^\perp \defas \M{I} - \MT{B} \, \M{W} \, \M{B} = \M{0}$
should be a matrix containing only zeros. However, this is not the
case, due to the loss of orthogonality in the three term
relationship resulting from numerical errors. These numerical
errors determine the quality of the basis functions and for which
we require a measure. The determinant of $\M{G}$ has in the past
been used as a measure for the quality $\epsilon_g = \mathrm{det}
\,\M{G}$ of the basis functions. However, this measure does not
yield stable estimates~\cite[Chapter 2, Sec. 2.7.3]{Golub}. We
propose the Frobenius norm of $\M{G}^\perp$ as an error measure,
i.e., $\epsilon_F = \| \M{G}^\perp \|_F$, this is the sum of the
square of all errors w.r.t. the orthogonality of the basis
functions, $\epsilon_F \geq 0$. This is a posteriori measure,
i.e., we compute the basis functions and then determine their
quality. Wilkinson~\cite{Wilkinson1971} points out that a-priory
prediction of error bounds yield unreliable results and a
posteriori analysis is preferred. The numerical results obtained
for different synthesis procedures can be found in
Section~(\ref{sec:quality}).
%
\subsection{Numerically Stable Synthesis of Basis Functions and
their Derivatives}
Gram~\cite{Gram1883} proposed what is now known as the
Gram-Schmidt orthogonalization process to generate
polynomials~\cite{Barnard1998}. The Gram-Schmidt process is,
however, numerically unstable~\cite[Chapter 5]{Golub} and errors
accumulate as the number of integrations increases, i.e., with
increasing polynomial degree. This precludes the synthesis of
polynomials of higher degree with this method. Considerable
research has been performed on discrete polynomials and their
synthesis~\cite{Mukundan2001,Yap2003,Yap2005,Yang2006,Hosny2007,Zhu2007,Zhu2007B,baik2007}.
The research was primarily in conjunction with the computation of
moments for image processing. None of these papers present a
method which is capable of synthesizing discrete orthogonal
polynomials of high quality for arbitrary nodes located within the
unit circle on the complex plane.

Here it is proposed to synthesize the polynomial basis functions
using a Lanczos process with complete
reorthogonalization~\cite[Chapter 9, p.
482]{Golub},\cite{oleary2008b}. The procedure can be summarized as
follows: Given a vector $\V{x}$ of $n$ nodes with mean $\bar{x}$,
i.e., the points at which the differential equation is to be
solved: first compute the two basis functions $\V{b}_0$, $\V{b}_1$
and initialize the matrix of basis functions $\M{B}$,
\begin{equation}
    \V{b}_0 = \V{1} / \sqrt{n}
    \hspace{1cm}
    \V{b}_1 = \frac{\V{x} - \bar{x}}{\norm{ \V{x} - \bar{x} }}
    \hspace{1cm}
    \text{and}
    \hspace{1cm}
    \M{B} = \left[\V{b}_0, \V{b}_1 \right].
\end{equation}

The remaining polynomials are synthesized by repeatedly performing
the following computations:
\begin{enumerate}
    \item Compute the polynomial of the next higher
    degree\footnote{The symbol $\circ$ represents the Hadamard product.},
    \begin{equation}
        \V{b}_n = \V{b}_1 \circ \V{b}_{n-1};
    \end{equation}
    \item perform a complete reorthogonalization,
    \begin{align}
        \V{b}_n &= \V{b}_n - \M{B} \, \MT{B} \, \V{b}_{n}\\
                &= \left\{ \M{I} - \M{B} \, \MT{B} \right\} \V{b}_{n}
    \end{align}
    by projection onto the orthogonal complement of all previously
    synthesized polynomials. It is important to note that the
    reorthogonalization is w.r.t. to the complete set of basis
    functions, not just the previous polynomial.
    \item Normalize the vector,
    \begin{equation}
        \V{b}_n = \frac{\V{b}_n}{\norm{ \V{b}_n }},
    \end{equation}
    \item and augment the matrix of basis functions,
    \begin{equation}
        \M{B} = \left[\M{B}, \V{b}_n \right].
    \end{equation}
\end{enumerate}
This procedure yields a set of orthonormal polynomials from a set
of arbitrary nodes located within the unit circle on the complex
plane. Although in~\cite{golubMeurant} the Lanczos process is used
to compute discrete orthogonal polynomials, the authors seem to
have overseen the possibility (necessity) of using complete
reorthogonalization at each step of the polynomial synthesis.

By taking the derivative of the recurrence relationship w.r.t.
$x$, we obtain the equations required to simultaneously synthesize
the differentials of the polynomials. This procedure appears
in~\cite{Kopriva2009} for the Legendre and Chebyshev polynomials.
Here the method is generalized to the synthesis of polynomials
from arbitrary nodes. With this, the synthesis procedure delivers
a set of orthonormal basis functions $\M{B}$ and their derivatives
$\dot{\M{B}}$.
%
\subsection{Weighted Basis Functions}\label{sec:W}

A set of discrete basis functions in matrix form, $\M{B}_w$ are
orthogonal with respect to a weighting matrix $\M{W}$ if,
\begin{equation} \label{eqn:orthWRTW}
    \MT{B}_w\M{W}\M{B}_w = \M{I}.
\end{equation}
In the case of a weighting function $w(x)$ the weighting matrix is
given by $\M{W} = \diag\left\{{w(x_1) \ldots w(x_n)}\right\}$.
Given a set of orthonormal basis functions $\M{B}$ and a positive
definite weighting matrix $\M{W}$, there exists a set of weighted
basis functions $\M{B}_w$, such that $\M{B}_w = \M{B} \, \M{U}$,
whereby $\M{U}$ is a full rank upper triangular matrix.
Substituting into Equation~(\ref{eqn:orthWRTW}) yields,
\begin{equation}
    \MT{U}\, \MT{B} \, \M{W} \, \M{B} \, \M{U} = \M{I}.
\end{equation}
Since $\M{U}$ is full rank, we may invert it to obtain,
\begin{equation}
    \MT{B} \, \M{W} \, \M{B}  = \MTI{U} \, \M{U}^{-1}.
\end{equation}
The Cholesky decomposition $\chol{\M{A}}$ of a matrix exists and
is unique such that $\M{A} = \M{G}\,\MT{G}$ if $\M{A}$ is real
positive definite. The matrix $\M{G}$ is a full rank lower
triangular matrix. Consequently, the Cholesky decomposition
$\chol{\MT{B} \, \M{W} \, \M{B}}$ exists if $\M{W}$ is real
positive definite, since $\M{B}$ is orthonormal. Applying the
decomposition yields,
\begin{equation}
    \MT{B} \, \M{W} \, \M{B}  = \M{G} \, \MT{G} = \MTI{U} \, \M{U}^{-1}.
\end{equation}
The sought matrix $\M{U}$ is clearly given by,
\begin{equation}
    \M{U} = \MTI{G}.
\end{equation}
With this the weighted basis functions are fully defined. The
condition number of the basis functions depends soley on the
condition number of the weighting matrix $\M{W}$. In the case
where the weighting matrix is derived from a weighting function
$w(x)$, the condition number is determined by the extreme values
of $w(x)$.
%
\subsection{Differentiating Matrices}\label{sec:D}
There are both
global~\cite{Kesan2003,Kurt2008,Keskyn2011,Sezer1996,Kopriva2009}
and
local~\cite{Courant1928,Savitzky1964,Courant1967,Podlubny2000,Podlubny2009}
approaches to computing discrete estimates for derivatives. Global
methods proposed in the past have used the known relationship
between the coefficients of a polynomial and the coefficients of
the derivative of the polynomial to compute a differentiating
matrix.

The computation of a differentiating matrix from polynomial bases
proceeds as follows: The spectrum of the signal $\V{y}$ with
respect to the basis functions $\M{B}$ is computed as,
\begin{equation}
    \V{s} = \M{B}^+ \, \V{y}.
\end{equation}
For example, $\M{B}$ may be the Vandermonde matrix; this is case
with Taylor methods~\cite{Kesan2003,Kurt2008,Keskyn2011}. The
relationship between the spectrum $\V{s}$ and the spectrum of the
derivatives is given by,
\begin{equation}
    \dot{\V{s}} = \M{M} \, \V{s}
    \hspace{2cm}
    \text{whereby}
    \hspace{2cm}
    \M{M} =
    \begin{bmatrix}
        0 & 1 & 0 & \ldots & 0\\
        0 & 0 & 2 & \ldots & 0\\
        0 & 0 & 0 & \ddots & 0\\
        0 & 0 & 0 &\ldots & n\\
        0 & 0 & 0 &\ldots & 0
    \end{bmatrix}.
\end{equation}
Consequently,
\begin{equation}
    \dot{\V{y}} = \M{B} \, \M{M} \, \M{B}^+ \, \V{y} = \M{D} \, \V{y}.
\end{equation}
That is, the differentiating matrix is computed as, $\M{D} = \M{B}
\, \M{M} \, \M{B}^+$. In the case of the Taylor (Vandermonde)
matrix this involves computing the pseudo-inverse of the
Vandermonde matrix: with all the associated numerical problems. In
the case of the Chebyshev polynomials~\cite{Sezer1996,Kopriva2009}
$\M{B}^+ = \MT{B}$ and a different matrix $\M{M}$ is required,
see~\cite{Sezer1996} for details. The method is not appropriate if
arbitrary nodes are required, e.g. this may be required if the
framework is to be used to solve the problems associated with
monitoring mechanical structures~\cite{Oleary2012}. The advantage
of Global methods is that they deliver a differentiating matrix
which is valid for the complete support.

The solution chosen here is to compute $\M{D}$ form the basis
functions and their derivatives, i.e., given $\dot{\M{B}}$ and
$\M{B}$ an appropriate derivative operator, $\M{D}$, should have
the property that,
\begin{equation}
    \M{D} \, \M{B} = \dot{\M{B}}.
\end{equation}
Post-multiplying by $\MT{B}$ yields
\begin{equation}
    \M{D}_B \defas \M{D} \, \M{B} \, \MT{B} = \dot{\M{B}} \, \MT{B}.
\end{equation}
If the basis function set is complete, i.e., $\M{B} \, \MT{B}$
then the above equation yields the differentiating matrix
directly,
\begin{equation}\label{eqn:diffD}
    \M{D} = \dot{\M{B}} \, \MT{B}.
\end{equation}
This computation is valid for arbitrary nodes. If a truncated,
i.e., an incomplete, set of basis functions is used then $\M{B} \,
\MT{B}$ is the projection onto the the basis functions $\M{B}$ and
$\M{D}_B$ is then a regularizing differentiating matrix. The
matrix $\M{D}_B$ can be applied to the computation of estimates
for derivatives in the presence of noise.

\begin{figure}
    \centering
    \includegraphics[width=7cm]{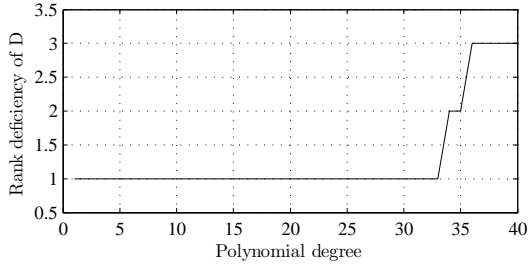}
    \caption{Rank deficiency of the differentiating matrix $\M{D}$ as a function of the
    degree of a Gram polynomial, when $\M{D}$ is synthesized using Equation~(\ref{eqn:diffD}).}
   \label{fig:RankDeficiencyD}
\end{figure}

A differentiating matrix should be rank-1 deficient; it should
have the constant vector as its null space, i.e., $\M{D} \, \V{1}
\, \alpha = \V{0}$. The properties of the $\M{D}$ are, however,
dependent on the nodes being used, e.g. the Chebyshev nodes permit
global differentiating matrices for very high degrees. With other
sets of nodes the condition number of a differentiating matrix can
increase with the degree of the polynomial being used. At some
point the matrix starts to have additional null spaces, which are
associated with numerical errors occuring due to insufficient
numerical precision, this effect is shown in
Figure~(\ref{fig:RankDeficiencyD}) for the Gram polynomials.

Once the condition number of a global differentiating matrix has
degenerated below an acceptable level, it becomes necessary to
compute local approximations~\cite{Savitzky1964,oleary2010C}.
Courant~\cite{Courant1928,Courant1967} proposed, in $1927$, using
both forward and backward differences to compute estimates for the
first derivative. The method has also been used
in~\cite{Podlubny2000,Podlubny2009} to this end. More commonly the
tri-diagonal matrix, shown here for $6$ points,
\begin{equation}
    \M{D}_t =
    \frac{1}{h}
    \begin{bmatrix}
   -1& 1& 0& 0& 0& 0\\
   -0.5& 0& 0.5& 0& 0& 0\\
   0& -0.5& 0& 0.5& 0& 0\\
   0& 0& -0.5& 0& 0.5& 0\\
   0& 0& 0& -0.5& 0& 0.5\\
   0& 0& 0& 0& -1& 1
    \end{bmatrix}
\end{equation}
is used to compute a  discrete local estimate for the first
differential\footnote{This is the matrix embedded in the Matlab
function \lstinline{gradient}.}. This operator is only of degree
$d=2$ accurate in the core of the approximation and at both ends
of the support only of degree $d=1$ accurate. This makes this
discrete operator unsuitable for the computation of derivatives at
the end of the support, as is required for BVPs and IVPs. It is
also not suitable for systems whose solutions are locally of
degree higher than $d > 2$. Furthermore, this operators assumes
equally spaced nodes. In~\cite{Amodio2011} higher order finite
difference schemes are proposed with end-point formulas. However,
only equidistant spaced nodes are considered. For example, the
appropriate three point operator for the Gram $\M{D}_{G,3}$ and
for the Chebyshev $\M{D}_{C,3}$ nodes, are,
\begin{equation}
    \M{D}_{G,3} =
    \begin{bmatrix}
    -4.5& 6& -1.5& 0& 0& 0\\
    -1.5& 0& 1.5& 0& 0& 0\\
    0& -1.5& 0& 1.5& 0& 0\\
    0& 0& -1.5& 0& 1.5& 0\\
    0& 0& 0& -1.5& 0& 1.5\\
    0& 0& 0& 1.5& -6& 4.5
    \end{bmatrix}
\end{equation}
and
\begin{equation}
    \M{D}_{C,3} =
    \begin{bmatrix}
    -5.2779& 6.0944& -0.8165& 0& 0& 0\\
    -2.4495& 1.633& 0.8165& 0& 0& 0\\
    0& -1.1954& 0.29886& 0.89658& 0& 0\\
    0& 0& -0.89658& -0.29886& 1.1954& 0\\
    0& 0& 0& -0.8165& -1.633& 2.4495\\
    0& 0& 0& 0.8165& -6.0944& 5.2779
    \end{bmatrix},
\end{equation}
both computed for $n = 6$ points in the interval $-1 < x < 1$.
Note the three points formulas at both ends of the support.

In keeping with the formulation of the basis functions for
arbitrary nodes: the method for local differential approximation
is also formulated here for arbitrary nodes. A generalized
formulation of local differentiating matrix requires the vector
$\V{x}$ of $n$ arbitrarily placed nodes, the support length $l_s$
and the degree $d$ of the approximation. Only odd support lengths
are considered here, to avoid the need for forward and backward
formulas. It is convenient to define the support length $l_s = 2
w_s + 1$ in terms of the half-width $w_s$. The vector $\V{x}$ of
nodes is segmented into $m = n - 2 w_s$ overlapping segments, for
each segment,
\begin{equation}
    \V{s}(i) = \V{x}(i - w_s : i + w_s)
    \hspace{1cm}
    \forall \, i \in [w_s + 1, n - w_s]
\end{equation}
a local set of basis functions $\M{B}_s$ and derivatives of the
basis functions $\dot{\M{B}}_s$ are computed. Then the
differentiating matrix associated with the segment is determined
$\M{D}_s = \dot{\M{B}}_s \, \MT{B}_s$. The first and last segments
yield the end-point formulas as required. The remaining segment
yields the required central formula of coefficients to locally
approximate the derivative. Clearly, for the inner-segments it is
only necessary to compute the center row vector of the local
differentiating operator $\M{D}_s$. The use of approximating or
interpolating polynomials leads to the generation of
differentiating matrices with and without regularization
respectively.
%
%
The Wilkinson diagram for the general structure of a local
differentiating matrix $\M{D}_L$ is shown in
Equation~(\ref{eqn:DMatStructure}) for the example of $l_s = 5$
and $n = 10$. The specific entries in the matrix are a function of
the spacing of the nodes.
\begin{equation}\label{eqn:DMatStructure}
    \M{D}_{5} =
    \begin{bmatrix}
    \times & \times & \times & \times & \times & 0 & 0 & 0 & 0 & 0\\
    \times & \times & \times & \times & \times & 0 & 0 & 0 & 0 & 0\\ \hline
    \times & \times & \times & \times & \times & 0 & 0 & 0 & 0 & 0\\
    0 & \times & \times & \times & \times & \times & 0 & 0 & 0 & 0\\
    0 & 0 & \times & \times & \times & \times & \times & 0 & 0 & 0\\
    0 & 0 & 0 & \times & \times & \times & \times & \times & 0 & 0\\
    0 & 0 & 0 & 0 & \times & \times & \times & \times & \times & 0\\
    0 & 0 & 0 & 0 & 0 & \times & \times & \times & \times & \times\\ \hline
    0 & 0 & 0 & 0 & 0 & \times & \times & \times & \times & \times\\
    0 & 0 & 0 & 0 & 0 & \times & \times & \times & \times & \times
    \end{bmatrix}
\end{equation}
All computations of the local derivative are of length $l_s$ and
of constant approximation degree $d_a = 2 w_s$ over the complete
support. This is important if derivatives are to be computed at
the ends of the support; furthermore, errors at the end of the
support associated with inconsistent approximations will propagate
through the entire solution when $\M{D}$ is being used in the
solution of differential equations. This procedure proposed here
delivers a local differentiating matrix for arbitrary nodes.
\subsection{Defining Constraints}\label{sec:C}
In Section~(\ref{sec:D}) it was shown that a discrete
approximation to differentiation can be computed as a linear
matrix operator. Consequently, both differential and integral
constraints are linear. In the framework proposed here, a
constraint is implemented by restricting a linear combination
$\VT{c} \, \V{y}$ of the solution vector $\V{y}$ to have a scalar
value $d$, i.e.,
\begin{equation}
    \VT{c} \, \V{y} = d.
\end{equation}
This is a very general mechanism, since any constraining function
can be implemented at a point $x_i$ for which a linear $n$ point
expansion around this point exists.
%
%
To give an example, consider the $C^2$ continuous periodicity
constraint $y(0) = y(1)$, $\dot{y}(0) = \dot{y}(1)$ and
$\ddot{y}(0) = \ddot{y}(1)$: given the differentiating matrix
$\M{D}$ and $\M{D}^2$, the three constraints can be formulated as:
\begin{align}
    [1, 0, \ldots , 0, -1] \, \V{y} = \VT{c}_1 \, \V{y} &= 0,\\
    \left\{\M{D}(1,:) - \M{D}(end,:)\right\} \, \V{y} = \VT{c}_2 \, \V{y} &= 0,\\
    \left\{\M{D}^2(1,:) - \M{D}^2(end,:)\right\} \, \V{y} = \VT{c}_3 \, \V{y} &=
    0.
\end{align}
Given a set of $m$ constraints, the constraining vectors $\V{c}_i$
are concatenated to form the matrix $\M{C} = [\V{c}_1 \ldots
\V{c}_m]$ and the corresponding scalars $d_i$ form the vector
$\VT{d} = [d_1 \ldots d_m]$, so that,
\begin{equation}
    \M{C} \, \V{y} = \V{d}.
\end{equation}

\subsection{Homogeneously Constrained Admissible Functions}
\label{sec:HCAF}

Starting from a set of basis functions $\M{B}$ such that $\MT{B}
\, \M{W} \, \M{B} = \M{I}$, we wish to derive a method of
synthesizing a set of constrained basis functions $\M{B}_c$ which
fulfil the conditions:
\begin{equation}\label{eqn:conditions}
    \MT{B}_c \, \M{W} \, \M{B}_c = \M{I},
    \hspace{1cm}
    \MT{C} \, \M{B}_c = \V{0}
    \hspace{1cm}
    \text{and}
    \hspace{1cm}
    \M{B}_c = \M{B} \, \M{X},
\end{equation}
i.e., the constrained basis functions form an orthonormal basis
set with respect to the weighting matrix $\M{W}$. If $\M{B}$ is
orthonormal, i.e., $\MT{B} \, \M{B} = \M{I}$ then so is $\M{B}_c$.
The constrained basis functions fulfil the homogeneous constraints
defined by $\M{C}$. If $\M{B}$ is complete then it spans the
complete $n \times n$ space, given $p = \rank{(\M{C})}$, i.e., the
number of independent constraints, $\M{B}_c$ is of dimension $ n
\times (n - p)$ and spans the complete space in which the
constraints are fulfilled. Consequently, all possible vectors
$\V{y}$ which fulfil the constraints are given by,
\begin{equation}\label{eqn:allYs}
    \V{y} = \M{B}_c \, \V{\alpha}
\end{equation}
where $\V{\alpha}$ is an $n - p$ vector.

A solution to the task of determining $\M{X}$ was presented
in~\cite{Oleary2012}; however, a more succinct derivation is
provided here. The conditions from Equation~(\ref{eqn:conditions})
require,
\begin{equation}
    \MT{C} \, \M{B} \, \M{X} = \V{0}
\end{equation}
and with this $\M{X}$ must lie in the null space of $\MT{C} \,
\M{B}$. Applying QR decomposition to $\MT{B} \, \M{C}$ yields,
\begin{equation}
    \M{Q} \, \M{R} = \MT{B} \, \M{C},
\end{equation}
and consequently,
\begin{equation}
    \MT{X} \, \M{Q} \, \M{R} = \V{0}
\end{equation}
The matrices $\M{Q}$ and $\M{R}$ are partitioned according to the
span and null space of $\MT{B} \, \M{C}$,
\begin{equation}
    \M{Q} = [\M{Q}_s, \M{Q}_n]
    \hspace{1cm}
    \text{and}
    \hspace{1cm}
    \M{R} =
    \begin{bmatrix}
        \M{R}_s\\
        \M{0}
    \end{bmatrix}.,
\end{equation}
with $\M{R}_s$ of dimension $p \times p$. The $n \times p$ matrix
$\M{Q}_s$ forms a basis set for the span and the $n \times (n -
p)$ matrix $\M{Q}_n$ forms a basis set for the null space of
$\MT{B} \, \M{C}$. Consequently,
\begin{equation}\label{eqn:qr1}
    \MT{X} \, \M{Q}_s = \M{0}
    \hspace{1cm}
    \text{and}
    \hspace{1cm}
    \MTb{\MT{X} \, \M{Q}_n} \M{W} \, \MT{X} \, \M{Q}_n = \M{I}.
\end{equation}
Now applying an RQ decomposition to $\M{Q}_n$ yields,
\begin{equation}
    \hat{\M{R}} \, \hat{\M{Q}}_n = \M{Q}_n.
\end{equation}
$\hat{\M{R}}$ is orthonormal, since both $\hat{\M{Q}}_n$ and
$\M{Q}_n$ are by definition orthonormal.  Now, selecting $\M{X} =
\hat{\M{R}}$ yields $\MT{X} \, \hat{\M{R}} \, \hat{\M{Q}}_n =
\hat{\M{Q}}_n$, and with this all the conditions from
Equation~(\ref{eqn:conditions}) are fulfilled. The matrix $\M{X}$
being orthonormal ensures that $\M{B}_c$ fulfils the same
orthonormal condition as does $\M{B}$. Furthermore, $\M{X}$ has an
implicit partitioning,
\begin{equation}
    \M{X}
    =
    \begin{bmatrix}
        \M{X}_1\\
        \M{X}_2
    \end{bmatrix}
\end{equation}
whereby, $\M{X}_1$ is a $p \times (n - p)$ block matrix and
$\M{X}_2$ is a $(n - p) \times (n - p)$ upper triangular matrix.
This structure ensures that the number of roots in the constrained
basis functions $\M{B}_c$ is ordered in the same manner as in
$\M{B}$.

\section{Discretizing and Solving Ordinary Differential
Equations}\label{sec:dsolve}
In the previous section all the matrices required for the
discretization of ordinary differential equations were derived. In
this section the discretization of initial-, boundary- and inner
value problems is presented together with the associated methods
of solving the resulting matrix equations.
\subsection{Initial Value Problems}
In this paper we are considering the solution of linear ordinary
differential equations with constant or variable coefficients,
they can in general be formulated as,
\begin{equation}\label{eqn:deDef}
    p_k(x) \, y^{(k)}(x) \ldots + p_1(x) \, y^{(1)}(x) + p_0(x) \, y(x) =
    g(x)
\end{equation}
to which a set of $k$ constraints are required to ensure a unique
solution. The term $y^{(k)}(x)$ represents the $k^{\text{th}}$
derivative of $y(x)$. Given the matrices derived previously, the
discretization of Equation~(\ref{eqn:deDef}) is direct and simple,
each term $p_k(x) \, y^{(k)}(x)$ is dicreteized as follows: The
matrix $\M{P}_k$ is formed such that $\M{P}_k =
\diag{\left\{p_k(\V{x})\right\}}$, whereby $p_k(\V{x})$ is the
vector of values obtained by evaluating the function $p_n(x)$ at
the vector of points $\V{x}$; the term $y^{(k)}(x)$ is discretized
as $\M{D}^k \, \V{y}$, i.e., the $k^{\text{th}}$ power of $\M{D}$,
which is the differentiating matrix derived in
Section~(\ref{sec:D}). Summarizing, each term is discretized as
follows,
\begin{equation}
    p_k(x) \, y^{(k)}(x) \rightarrow \M{P}_k \, \M{D}^k \, \V{y}.
\end{equation}
and the vector $\V{g} = g(\V{x})$. Applying this to all terms in
Equation~(\ref{eqn:deDef}) yields,
\begin{equation}
    \M{P}_k \, \M{D}^k \, \V{y} \ldots
    + \M{P}_1 \, \M{D} \, \V{y}
    + \M{P}_0 \,  \V{y}
    =
    \V{g}
\end{equation}
The matrix equivalent of the linear  differential operator $\M{L}$
is now defined as,
\begin{equation}
    \M{L} \defas \M{P}_k \, \M{D}^k \ldots
    + \M{P}_1 \, \M{D}
    + \M{P}_0,
\end{equation}
and the set of $k$ constraints are implemented as defined in
Section~(\ref{sec:C}), yielding
\begin{equation}\label{eqn:sys2solve}
    \M{L} \,  \V{y} = \V{g}
    \hspace{1cm}
    \text{given}
    \hspace{1cm}
    \MT{C} \, \V{y} = \V{d}.
\end{equation}
the matrix $\M{C}$ has the dimension $n \times k$.

A unique solution to the ODE exists only if
\begin{equation}
    \rank{
    \begin{bmatrix}
        \M{L} \\
        \MT{C}
    \end{bmatrix}} = n
\end{equation}
i.e., the linear differential operator and the constraints must
form a full rank system of equations. There are many Engineering
application where this is not the case, e.g. the equations for the
vibration of a beam, and Sturm-Liouville problems. A different
solution approach is proposed for this class of problems in
Section~(\ref{sec:BVPSL}).
%
\subsubsection{Solution as a constrained least squares problem} The
formulation of determining $\V{y}$ from
Equation~(\ref{eqn:sys2solve}) as the solution of a least squares
minimization problem yields,
\begin{equation}\label{eqn:LS}
    \min_{y} \| \M{L} \, \V{y} - \V{g}\|_2^2
    \hspace{1cm}
    \text{given}
    \hspace{1cm}
    \MT{C} \, \V{y} = \V{d}.
\end{equation}
This is the well known problem of least squares with equality
constraints (LSE). Efficient and accurate solutions can be found
in~\cite[Chapter 12]{Golub}. This method will yield solutions for
ODEs with consistent constraints and a least squares solution in
the case of over-constrained systems and perturbed systems. It is
not a suitable approach for Sturm-Liouville type problems.

The worst case number of floating point operations (FLOPS)
required to perform the computation is known a-priori. This, by
definition, makes the method suitable for real time applications.
\subsubsection{Spectral Reqularization}

Spectral regularization is introduced here to limit the number of
zeros in the basis functions and with this to reduce the errors
associated with aliasing. Assuming $\V{y}$ can be sufficiently
accurately approximated by a series of $r$ orthonormal basis
functions, we may write,
\begin{equation}
    \V{y} = \M{B}_r \, \V{\alpha},
\end{equation}
whereby $\M{B}_r = \M{B}(:,1\ldots r)$. Now defining $\M{L}_r
\defas \M{L} \, \M{B}_r$ and $\M{C}_r \defas \MT{B}_r \, \M{C}$,
and substituting into Equation~(\ref{eqn:LS}) yields,
\begin{equation}\label{eqn:regsys2solve}
    \min_{\V{\alpha}} \| \M{L}_r \, \V{\alpha} - \V{g}\|_2^2
    \hspace{1cm}
    \text{given}
    \hspace{1cm}
    \MT{C}_r \, \V{\alpha} = \V{d},
\end{equation}
whereby the series coefficients $\V{\alpha}$ are to be determined.
In addition to introducing regularization, the truncated basis
functions also reduce the size of the LS problem to be solved.
\subsubsection{Solution of Homogeneously Constrained IVPs}
Homogeneously Constrained IVPs for a special subclass of problems
for which there is a particularly simple solution. Let the
solution $\V{y}$ be a linear combination of a set of constrained
basis functions, i.e., $\V{y} = \M{B}_c \, \V{\alpha}$, which
fulfil the homogeneous constraints $\MT{C} \, \M{B}_c = 0$
associated with the IVP. Equation~(\ref{eqn:LS}) now simplifies to
the unconstrained least squares problem,
\begin{equation}
    \min_{\V{\alpha}} \| \M{L} \, \M{B}_c \, \V{\alpha} -
    \V{g}\|_2^2.
\end{equation}
The solution of which is,
\begin{equation}
    \V{\alpha} =
    \left\{
        \M{L} \, \M{B}_c
    \right\}^+ \, \V{g}
\end{equation}
since $\mathrm{null} \left\{\M{L} \, \M{B}_c \right\} = \V{0}$ if
a unique solution exists. Consequently,
\begin{equation}
    \V{y} = \M{B}_c \,
    \left\{
        \M{L} \, \M{B}_c
    \right\}^+ \, \V{g}
\end{equation}
\subsection{Sturm-Liouville and Boundary Value Problems}\label{sec:BVPSL}
A Sturm-Liouville problem is a second order ODE with the following
structure,
\begin{equation}\label{eqn:StLio}
    -\frac{\md}{\md x} \,
    \left[
        p(x) \, \frac{\md y}{\md x}
    \right]
     + g(x) \, y
    = \lambda \, w(x) \, y,
\end{equation}
in the finite interval $x_1 \leq x \leq x_n$, where $p(x)$, $g(x)$
and $w(x)$ are real-valued strictly positive. Additionally there
are two boundary conditions which are most commonly formulated as,
\begin{align}
    a_1 \, y(x_1) + a_2 \, \dot{y}(x_1) &= 0,\\
    b_1 \, y(x_2) + b_2 \, \dot{y}(x_2) &= 0.
\end{align}
There are some important properties of Sturm-Liouville
equations~\cite{Ledoux2007} which must be considered when
implementing a discrete solution:
\begin{enumerate}
    \item All eigenvalues are real and there is no largest
    eigenvalue, i.e., there are an infinite number of eigenvalues
    and $\lambda_m \rightarrow \infty$ as $m \rightarrow \infty$.
    Given a set of $n$ discrete points $\V{x}$ there can
    theoretically only be $n$ eigenvalues;
    \item The $m^{\text{th}}$ eigenfunction has $m$ zeros on
    the interval $a < x < b$. However, given $n$ points (samples)
    only functions with a maximum of $n/2$ zeros can be discribed
    without aliasing. Consider the
    Sturm-Liouville equation $\ddot{y} - \lambda \, {y} = 0$ with
    the constraints
    $y(0)=0$ and $y(\pi) = 0$. This equation is known to have the
    eigenfunctions $\Phi_m(x) = \sqrt{2}\,\sin{(m \, \pi \, x)}$.
    Consequently, a discrete solution can only model the first $n/2$ eigenpairs
    correctly.

%
    \item The eigenfunctions are orthogonal with respect to the
    weighting function $w(x)$, i.e., $\int_a^b w(x) \, \Phi_i(x)
    \Phi_j(x)= \delta(i,j)$.
\end{enumerate}

The general Sturm-Liouville problem formulated in
Equation~(\ref{eqn:StLio}) with its corresponding boundary
conditions can be discritized directly as,
\begin{equation}\label{eqn:DStLio}
    \left\{
        \M{D}\,\M{P}\,\M{D} - \M{G}
    \right\} \, \V{y}
    = - \lambda \, \M{W} \, \V{y}
    \hspace{1cm}
    \text{given}
    \hspace{1cm}
    \MT{C} \, \V{y} = \V{0}.
\end{equation}
whereby, $\M{P} = \diag{\left\{p(\V{x})\right\}}$, $\M{G} =
\diag{\left\{g(\V{x})\right\}}$ and $\M{W} =
\diag{\left\{w(\V{x})\right\}}$. A direct solution of this
equation will, however, yield unstable results due to aliasing.

We now introduce a set of weighted and constrained basis functions
$\M{B}_w$ which fulfil the orthogonality condition $\MT{B}_w \,
\M{W} \, \M{B}_w = \M{I}$ and boundary conditions $\MT{C} \, \M{B}
= \V{0}$. These basis functions are admissible functions for the
Sturm-Liouville problem. The number of zeros in the basis
functions increases from left to right in the matrix. The number
of zeros in the admissible functions is limited, so as to avoid
aliasing, by truncating to the first $k = n/2$ basis functions,
i.e., $\M{B}_a = \M{B}_w(:,1:k)$. The eigenfunctions are now found
as linear combinations of these admissible functions, i.e., $\V{y}
= \M{B}_a \, \V{\alpha}$. Substituting this into
Equation~(\ref{eqn:DStLio}) yields,
\begin{equation}
    \left\{
        \M{D}\,\M{P}\,\M{D} - \M{G}
    \right\} \, \M{B}_a \, \V{\alpha}
    = - \lambda \, \M{W} \, \M{B}_a \, \V{\alpha}.
\end{equation}
Pre-multiplying both sides by $\MT{B}_a$ now yields,
\begin{equation}
    \MT{B}_a \, \left\{
        \M{D}\,\M{P}\,\M{D} - \M{G}
    \right\} \, \M{B}_a \, \V{\alpha}
    = - \lambda \, \V{\alpha}.
\end{equation}
since $\MT{B}_a \, \M{W} \, \M{B}_a = \M{I}$. Now defining
$\M{L}_a
\defas \MT{B}_a \, \left\{\M{D}\,\M{P}\,\M{D} + \M{G} \right\} \,
\M{B}_a$ yields a standard eigenvector problem,
\begin{align}
    \left\{\M{L}_a \, + \lambda \, \M{I}\right\} \, \V{\alpha} &=
    0\label{eqn:eigSol1},\\
    \V{y} &= \M{B}_a \, \V{\alpha}\label{eqn:eigSol2}.
\end{align}
Solving Equation~(\ref{eqn:eigSol1}) for the eigenvalues
$\lambda_i$ and the eigenvectors $\V{\alpha}_i$, then back
substituting $\V{\alpha}_i$ into Equation~(\ref{eqn:eigSol2})
yields the desired eigenfunctions.

It is important and interesting to note the the matrix $\M{L}_a$
is of dimension $n/2 \, \times \, n/2$, in contrast to the
original matrix $\M{L} = \M{D}\,\M{P}\,\M{D} - \M{G}$ which is of
dimension $n \times n$. Consequently, dealing with the aliasing
has also reduced the size of the eigenvalue problem to be solved.
In the worst case an eigen-decomposition is of
complexity\footnote{Indeed there are more efficient algorithms;
however, their complexity depends on the structure of the matrix
and the distance between the eigenvalues. Consequently, no general
statements can be made about these methods.} between
$\mathcal{O}(n^2)$ and $\mathcal{O}(n^3)$. The improvement in
speed is then in the range of a factor of $4$ to $8$, while
simultaneously improving the accuracy of the solution. However,
some of the computation gains are spent on additional pre- and
post-calculations. A consequence of Equation~(\ref{eqn:eigSol2})
is that the matrix of eigenvectors $\M{\alpha}$, contains the
spectrum of the eigenfunctions with respect to the basis functions
used, i.e., the Rayleigh-Ritz coefficients.
\section{Performance Testing}\label{sec:test}
In this section a selection of examples are presented to
demonstrate the functionality of the proposed methods\footnote{A
MATLAB toolbox DOPbox is available at\\ 
http://www.mathworks.de/matlabcentral/fileexchange/41250 which
implements all the methods proposed in this paper; furthermore,
the code for all the following examples in also provided in the
toolbox: ODEbox
http://www.mathworks.de/matlabcentral/fileexchange/41354}.
\subsection{Quality of Basis Functions}\label{sec:quality}
\begin{figure}[h]
    \centering
    \subfigure[Complete polynomial basis sets, i.e., $d = n - 1$ where $n$ is the number of nodes.]{
    \includegraphics[width=6.1cm]{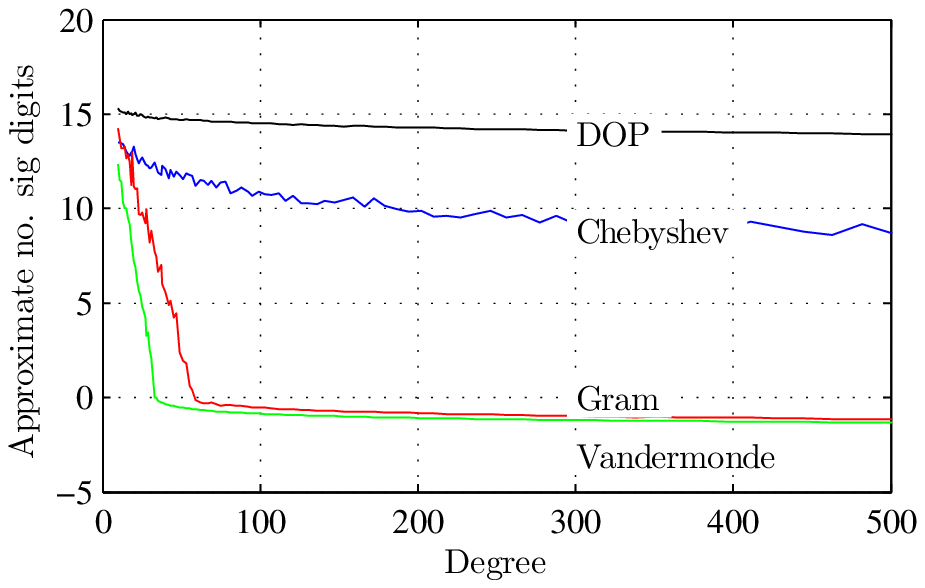}\label{fig:QualitySub1}}
    \hspace{3mm}
    \subfigure[Number of nodes $n = 1000$, the degree of the polynomial is increasing.]{
    \includegraphics[width=6.1cm]{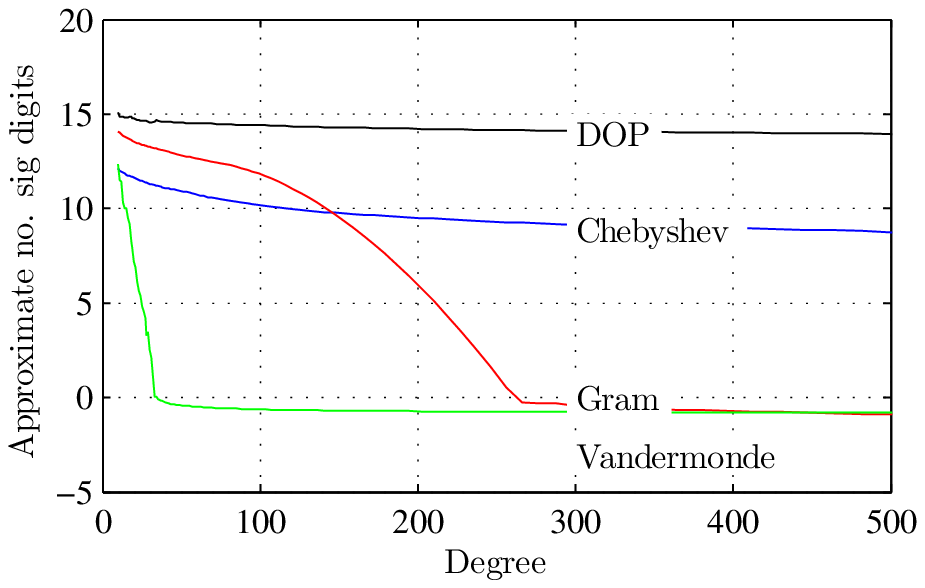}\label{fig:QualitySub2}}
    \caption{Estimate of the number of significant digits for different polynomials
     as a function of degree.
    DOP refers to the discrete orthonormal polynomial synthesis as proposed in this paper.}
   \label{fig:Quality}
\end{figure}

The first test addresses the quality of basis functions, since
these form the basis for all subsequent calculations. The
following polynomials are compared: a set of Gram polynomials
generated using Gram-Schmidt orthogonalization~\cite{Gram1883}; a
set of Chebyshev polynomials generated using the recurrence
relationship~\cite{numericalRecipes}; a Vandermonde matrix and a
set of polynomials synthesized using the method proposed in this
paper. The Frobenius norm of the projection onto the orthogonal
complement of the Gram matrix is used as an estimate of the total
error. The number of significant digits is then estimated to be $d
= - \log_{10}(\epsilon_F)$. Two computations were performed:
Figure~(\ref{fig:QualitySub1}) shows the result for complete
polynomial sets, i.e., the degree $d = n - 1$ where $n$ is the
number of nodes; Figure~(\ref{fig:QualitySub2}) is for a fixed
number of nodes $n = 1000$ and the degree of the polynomial is
progressively increased.
The results shown in Figure~(\ref{fig:Quality}) indicate that the
algorithm presented in this paper generates the most stable sets
of polynomials. For this reason the algorithm is used for the
generation of all bases required in this paper.
\subsection{Initial Value Problems}
In each of the following examples the analytical solution is
compared with the numerical results commuted using the newly
proposed method and a \textbf{Runga-Kutta} solution with variable
step size\footnote{See the MATLAB documentation for details of the
\textsf{ODE45} solver used for this computation.}. The
\lstinline{ODE45} is used for comparison in all the following
initial value problems.

\subsubsection{IVP Example 1}
\begin{figure}[h]
    \centering
    \subfigure[Comparison of the analytical solution,
    the numerical solutions using the new method (with support length $l_s =
    13$, the nodes are placed $x = 5 z^2$ and $z$ has $n = 85$
    evenly spaced points in the interval $ 0 \leq z \leq 1$)
    and using a Runge-Kutta procedure.
    The nodes are shown below the curve.]{
    \includegraphics[width=6.1cm]{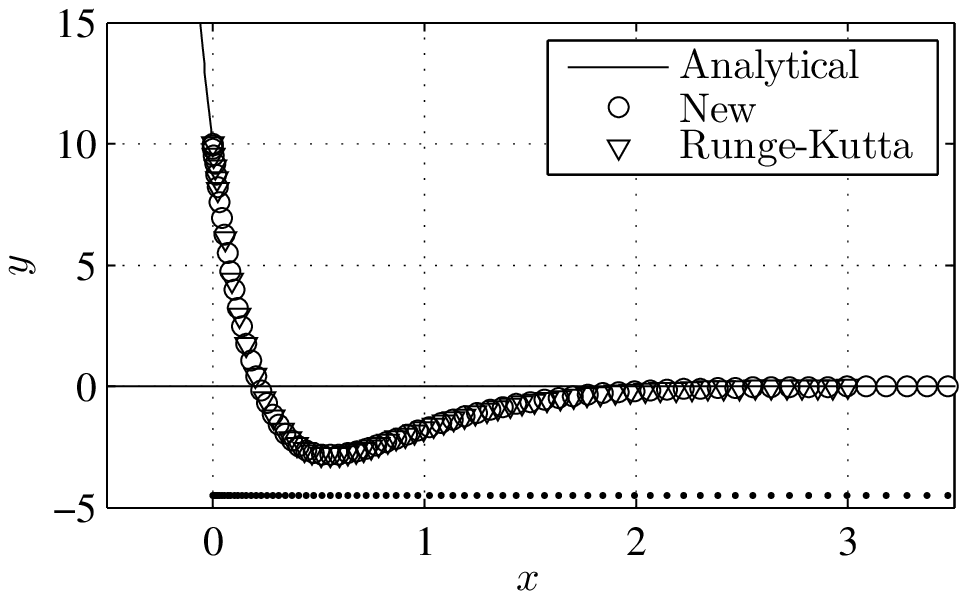}\label{fig:IVP1Sub1}}
    \hspace{3mm}
    \subfigure[Residual errors: (top) between the analytical solution
    and the new numerical procedure, (bottom) between the analytical solution
    and the numerical solution using the Runga-Kutta procedure. ]{
    \includegraphics[width=6.1cm]{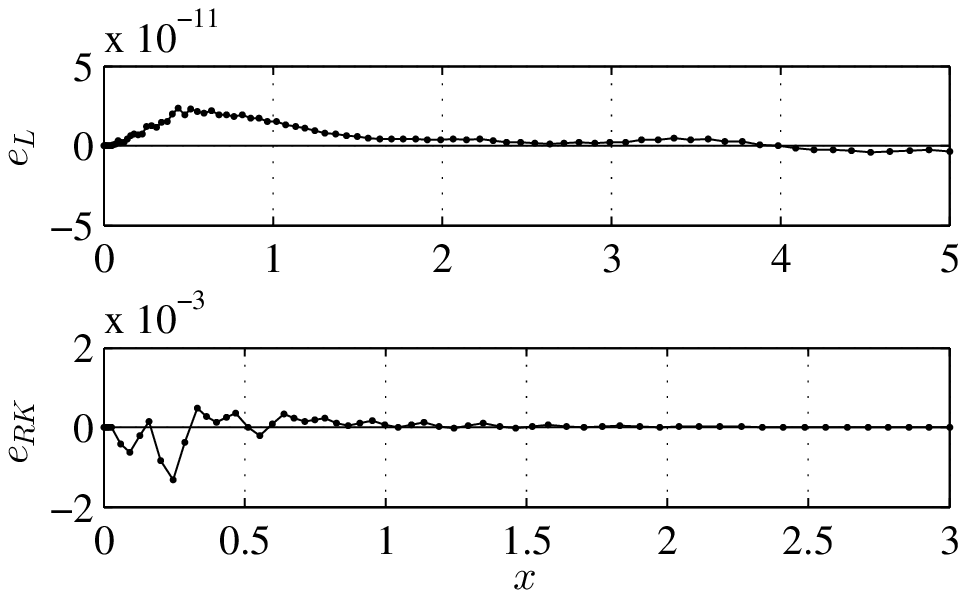}\label{fig:IVP1Sub2}}
    \caption{Computation results for the IVP Example 1 given in Equation~(\ref{eqn:IVP1}).}
    \label{fig:ivpExp1}
\end{figure}

The first example is a second order initial value problem with
constant coefficients, the equation is~\cite{Adams},
%
\begin{equation}\label{eqn:IVP1}
    \ddot{y} - 6\,\dot{y} - 9 y = 0,
    \hspace{5mm}
    \text{with},
    \hspace{5mm}
    y(0) = 10
    \hspace{5mm}
    \text{and}
    \hspace{5mm}
    \dot{y}(0) = -75.
\end{equation}
in the interval $0 \leq x \leq 3$ The analytical solution to this
equations is,
\begin{equation}
    y(x) = 10 \, e^{-3\,x} - 45 \,x \,{e^{-3\,x}}.
\end{equation}
The analytical solution and the results of the numerical solutions
are shown in Figure~(\ref{fig:IVP1Sub1}). Non-uniformly spaced
nodes were used, with a higher density of nodes where $\V{y}$ has
a higher first derivative. This demonstrates the possibility of
generating basis functions from arbitrary nodes. The residual
errors, see Figure~(\ref{fig:IVP1Sub2}), with the new method is
$7$ orders of magnitude smaller than with the \lstinline{ODE45}
method.
\subsubsection{IVP Example 2}
\begin{figure}[h]
    \centering
    \subfigure[Comparison of the analytical result,
    the numerical solutions using the new method (with support length $l_s =
    13$, there are $n = 73$ evenly spaced nodes in the interval $ 1 \leq x \leq 10$)
    and using a Runge-Kutta procedure.]{
    \includegraphics[width=6.1cm]{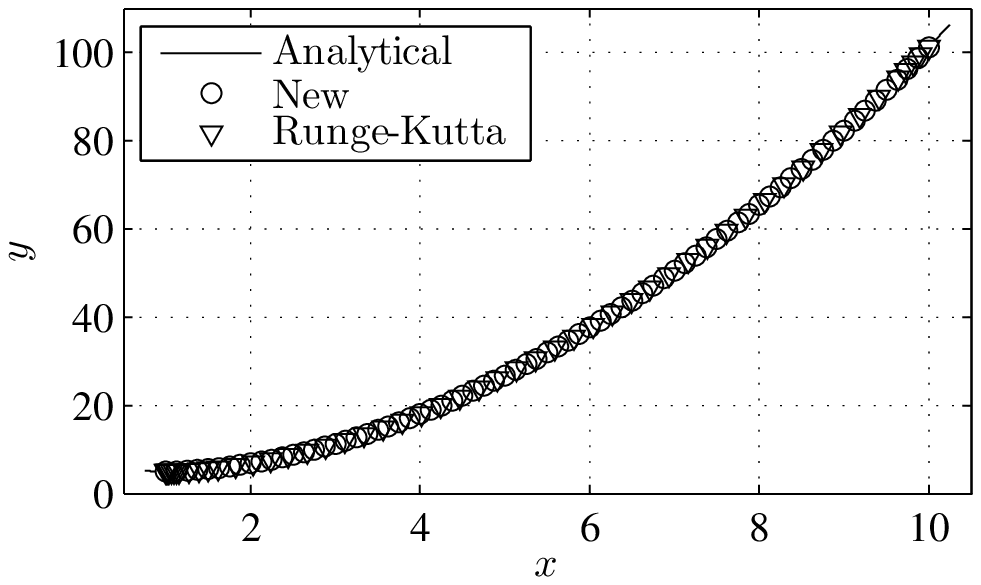}}
    \hspace{3mm}
    \subfigure[Residual errors: (top) between the analytical solution
    and the new numerical procedure, (bottom) between the analytical solution
    and the numerical solution using the Runga-Kutta procedure. ]{
    \includegraphics[width=6.1cm]{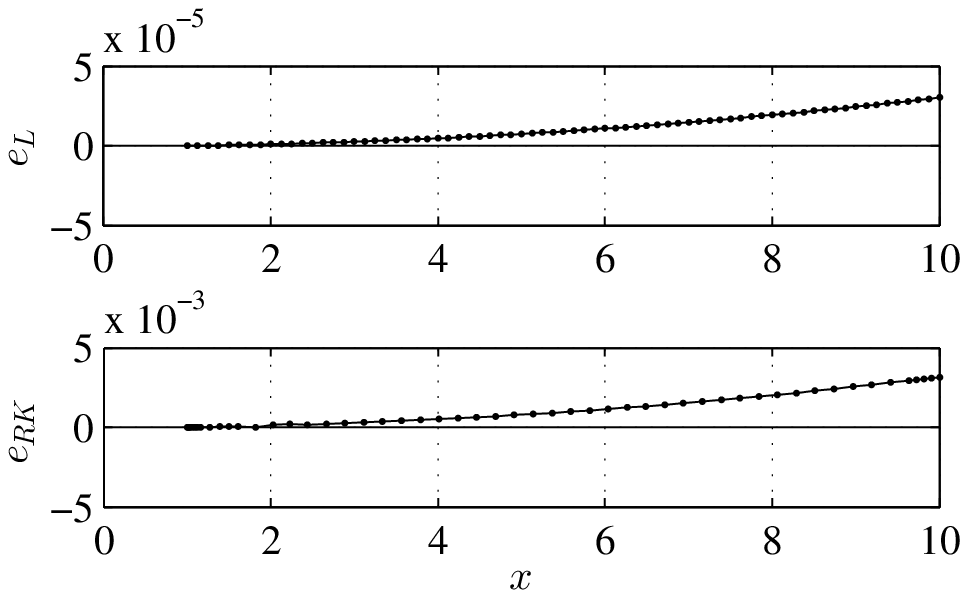}}
    \caption{Computation results for the IVP example 2 given in Equation~(\ref{eqn:IVP2}).}
    \label{fig:ivpExp2}
\end{figure}

The second example is a second order differential equation with
variable coefficients, the equation is~\cite{Adams},
\begin{equation}\label{eqn:IVP2}
    2\,x^2\,\ddot{y} - x\,\dot{y} - 2 y = 0,
    \hspace{5mm}
    \text{with},
    \hspace{5mm}
    y(1) = 5
    \hspace{5mm}
    \text{and}
    \hspace{5mm}
    \dot{y}(1) = 0.
\end{equation}
in the interval $1 \leq x \leq 10$. The analytical solution to
this equations is,
\begin{equation}
    y(x) = x^2 + \frac{4}{\sqrt{x}}.
\end{equation}
The comparison of the analytical solution with the numerical
solutions using the new method and a Runge-Kutta procedure is
shown in Figure~(\ref{fig:ivpExp2}). The residual error with the
new method is $3$ orders of magnitude smaller that with the
Runga-Kutta procedure. This example has demonstrated the ability
of the proposed method to solve initial value problems with
variable coefficients.
\subsubsection{IVP Example 3}
\begin{figure}[h]
    \centering
    \subfigure[Comparison of the analytical solution with,
    the numerical solutions using the new method (with support length $l_s =
    13$, there are $n = 73$ evenly spaced nodes in the interval $ 0 \leq x \leq 8$)
    and a Runge-Kutta procedure.]{
    \includegraphics[width=6.1cm]{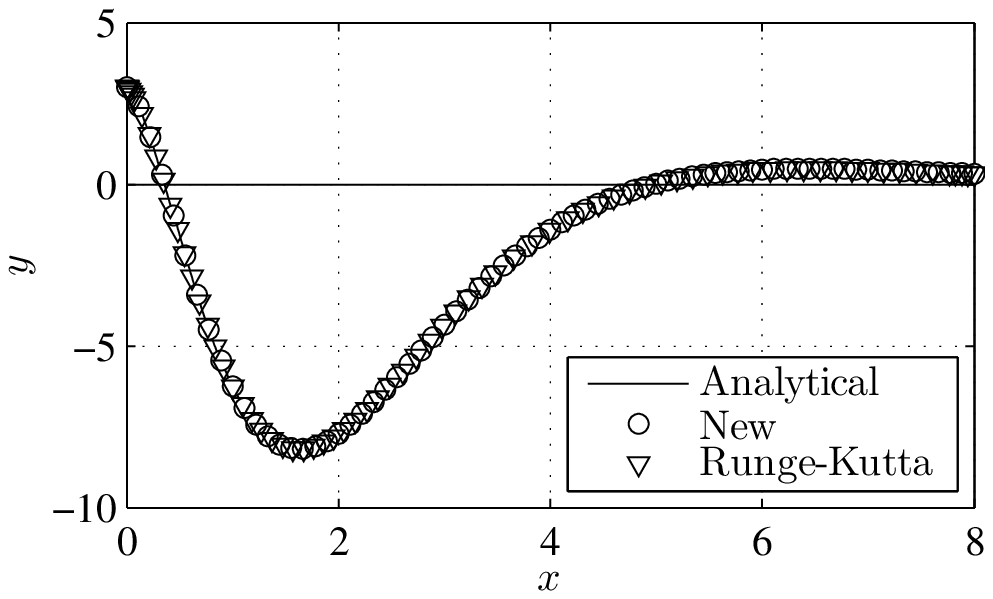}}
    \hspace{3mm}
    \subfigure[Residual errors: (top) between the analytical solution
    and the new numerical procedure, (bottom) between the analytical solution
    and the numerical solution using the Runga-Kutta procedure. ]{
    \includegraphics[width=6.1cm]{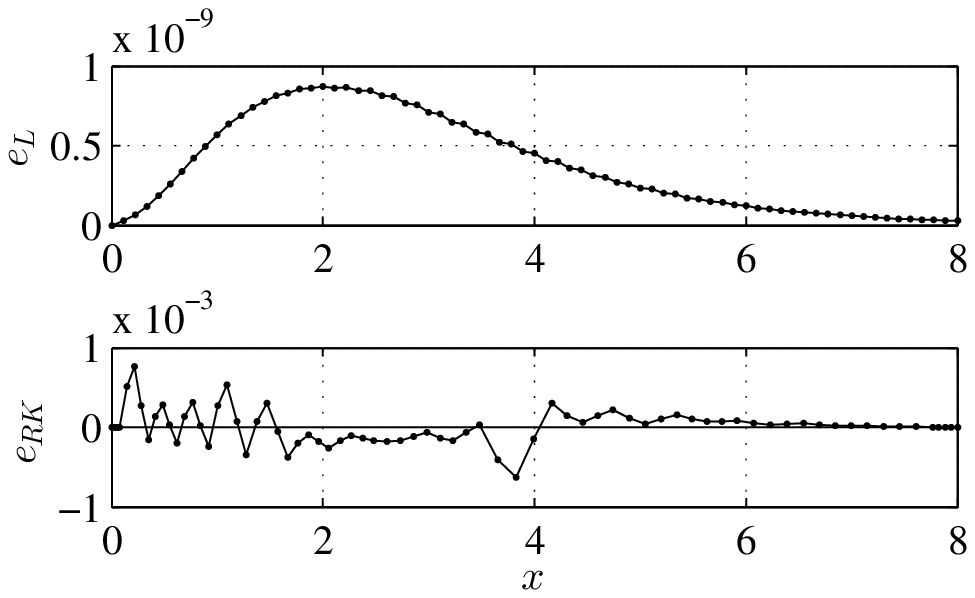}}
    \caption{Computation results for the IVP Example 2 given in Equation~(\ref{eqn:IVP3}).}
    \label{fig:ivpExp3}
\end{figure}

The third example~\cite{Adams} is a third order non-homogeneous
differential equation with constant coefficients, the equation is,
\begin{equation}\label{eqn:IVP3}
    \dddot{y} + 3\,\ddot{y} + 3 \dot{y} + y = 30 \, \me^{-x}
    \hspace{5mm}
    \text{with},
    \hspace{5mm}
    y(0) = 3,
    \hspace{5mm}
    \dot{y}(0) = -3,
    \hspace{5mm}
    \text{and}
    \hspace{5mm}
    \ddot{y}(0) = - 47
\end{equation}
in the interval $0 \leq x \leq 8$. The analytical solution to this
equation is,
\begin{equation}
    y(x) = (3 - 25\,x^2 + 5\,x^3)\,\me^{-x}.
\end{equation}
The comparison of the analytical solution with the numerical
solutions using the new method and a Runge-Kutta procedure is
shown in Figure~(\ref{fig:ivpExp3}). Once again the residual error
with the new method is orders of magnitude smaller that with the
Runga-Kutta procedure.

\subsection{Sturm-Liouville Problems}
The test package for Sturm-Liouville Solvers~\cite{Pryce1999} has
been used as a source of test cases in this section\footnote{At
this point we feel it is important to note that the framework
presented here is generally applicable to the solution of ODEs in
general and is not a dedicated Sturm-Liuville solver. The use of
Sturm-Liouville problems as a test cases is to demonstrate this
generality.}.
\subsubsection{Sturm-Liouville Example 1}
\begin{figure}[h]
    \centering
    \subfigure[Solution for $n = 100$.]{
    \includegraphics[width=6.1cm]{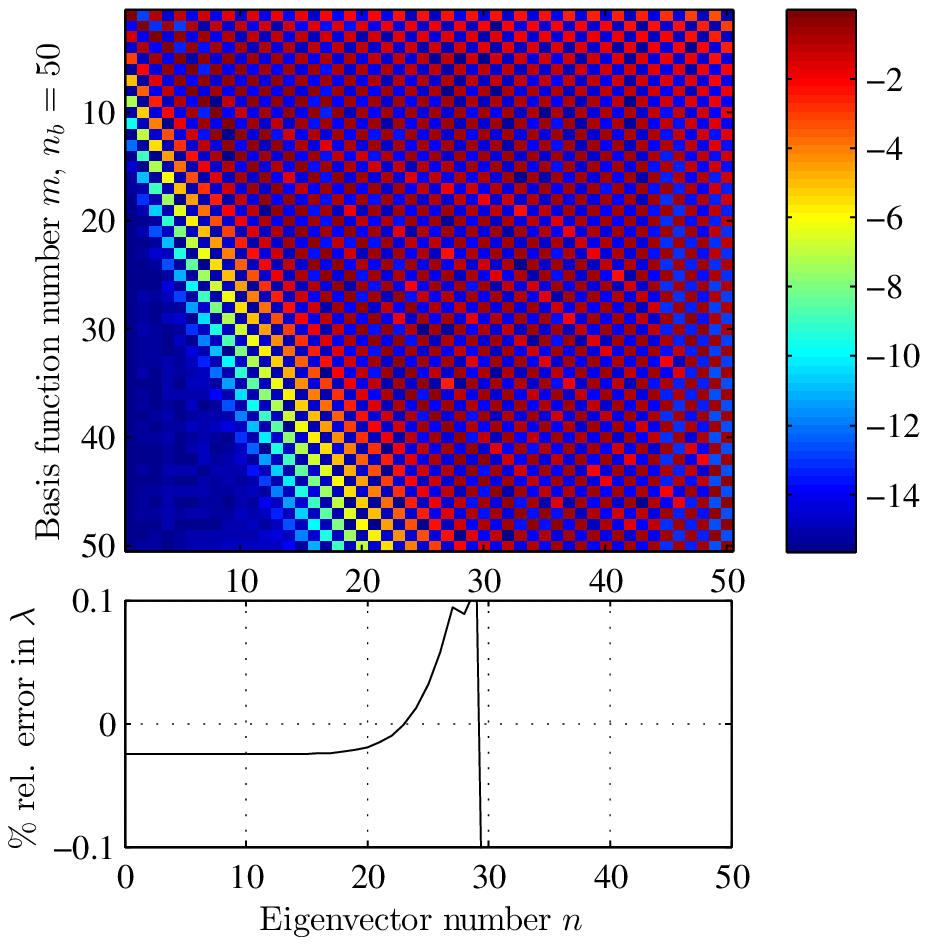}\label{fig:bvpExp1Sub1}}
    \hspace{3mm}
    \subfigure[Solution for $n = 1000$.]{
    \includegraphics[width=6.1cm]{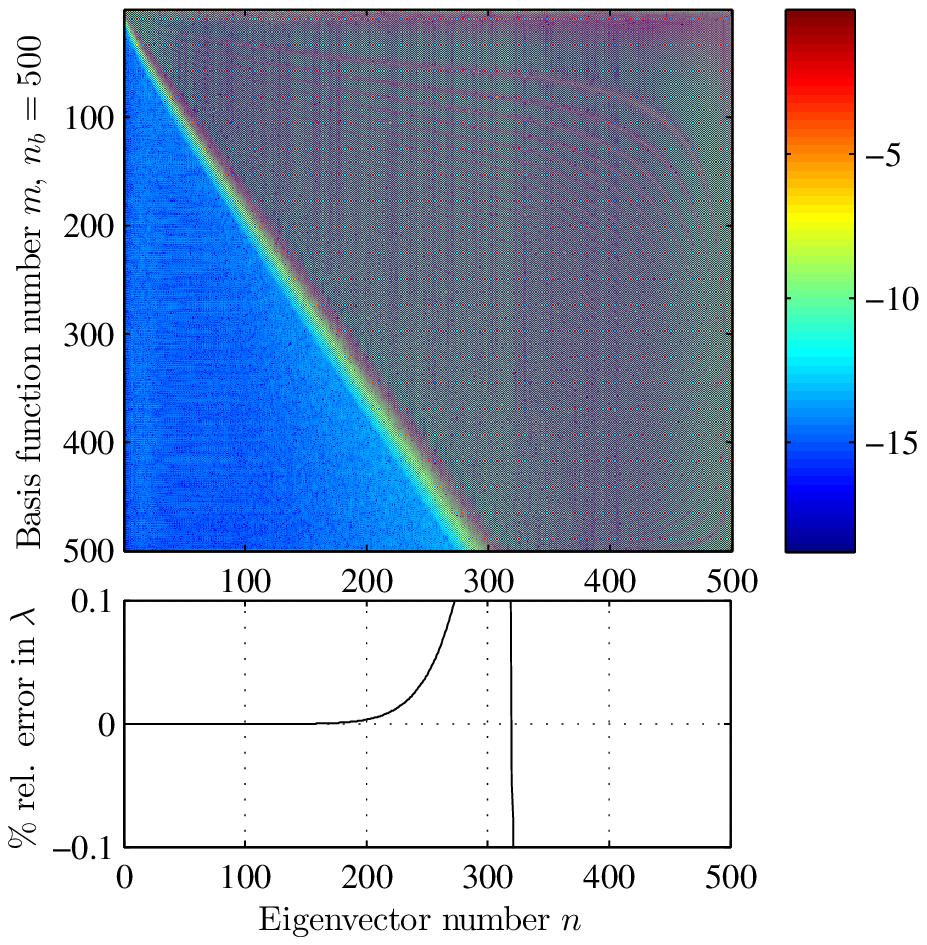}\label{fig:bvpExp1Sub2}}
    \caption{Solution of the Sturm-Liouville problem corresponding to the
    vibrating string for points. Top: spectrum of the eigenfunctions with respect to the admissible
    functions $(log_{10}(S)$. Bottom: The relative error in the eigenvalue $\lambda_k$ in $\%$.}
   \label{fig:bvpExp1}
\end{figure}

%
The simplest Sturm-Liouville problem is chosen as the first
example, since the analytical solution is known. This enables the
investigation of the stability of the numerical computation, i.e.,
how many eigenvalues can be computed to a given accuracy. It is
the equation of a vibrating string,
\begin{equation}
    \ddot{y} + \lambda y = 0,
    \hspace{5mm}
    \text{given}
    \hspace{5mm}
    y(0) = 0,
    \hspace{2mm}
    \text{and}
    \hspace{2mm}
    y(\pi) = 0
\end{equation}
in the interval $0 \leq x \leq \pi$. The analytical solution
yields the analytical eigenvalues $\lambda_k$ and eigenfunctions
$y_k$,
\begin{equation}
    \lambda_k = k,
    \hspace{5mm}
    y_k = \sin k x
    \hspace{5mm}
    \text{for}
    \hspace{5mm}
    k = 1\ldots \infty.
\end{equation}
The discrete solution has been computed on the interval $0 \leq x
\leq \pi$ sampled at the corresponding Chebyshev points; however,
scaled so that the first and last points lie exactly at $0$ and
$\pi$ respectively. For a given set of $n$ points $m = n/2$
constrained basis functions are computed which fulfil the boundary
values. These are the admissible functions used in what is
essentially a discrete equivalent of the Rayleigh-Ritz method. Two
different computations $n_1 = 100$ and $n_2 = 1000$, have been
performed to investigate the behavior of the solution with respect
to the number of points used. A support length $l_s = 13$ was used
during the generation of the differentiating matrix. The results
can be seen in Figure~(\ref{fig:bvpExp1Sub1})
and~(\ref{fig:bvpExp1Sub2}) respectively. The method returns the
coefficients of the series of admissible functions required to
generate each eigenfunction. The matrix of these coefficients is
denoted by $\M{R}$. We use the matrix $\M{S} =
\log_{10}\left\{\mathrm{abs}(\M{R}) \right\}$ as a visual
representation for the spectrum in the figures, since at
$\M{S}(i,j) \approx -16$ the numerical resolution of computation
environment is reached. This makes a visual recognition of when
the spectrum fails simple.


With $n_1 = 100$ approximately $k_1=28$ and for $n_2 = 1000$
approximately $k_2 = 280$ eigenvalues could be computed with a
relative error smaller than $0.1 \%$. This result is significantly
better than all previously reported results with matrix methods
for Sturm-liouville problems~\cite{Ledoux2007}. This confirms the
numerical stability of the approach. It also verifies that the
number of eigenfunctions which can be computed to a given accuracy
scales linearly with the number of nodes used.
\subsubsection{Sturm-Liouville Example 2}
\begin{figure}[h]
    \centering
    \subfigure[The first two eigenfunctions of the Mathieu differential equation, for $r= - 25$; note the negative sign for the value of $r$.]{
    \includegraphics[width=6.1cm]{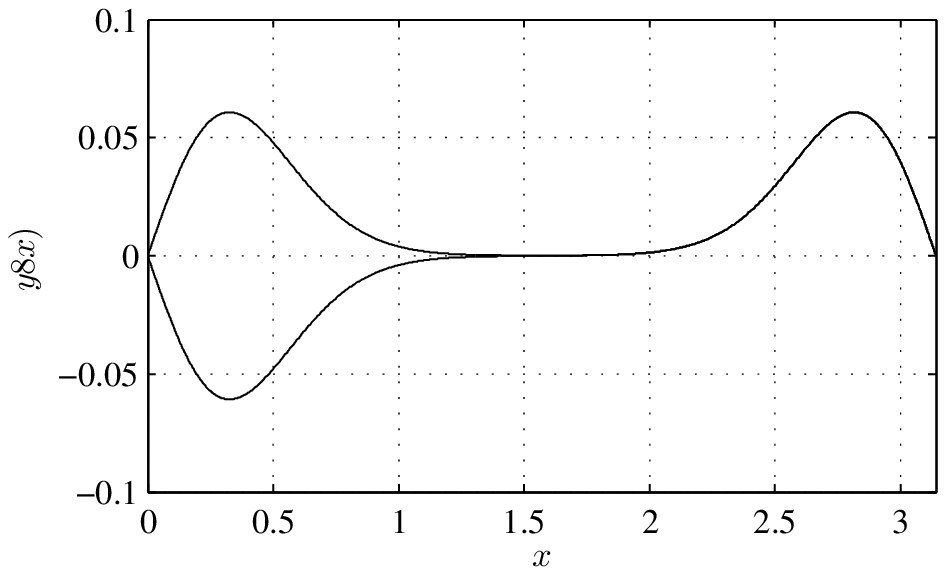}\label{fig:bvpExp2Sub1}}
    \hspace{3mm}
    \subfigure[Top: Rayleigh-Ritz Spectrum of the Mathieu differential equation. Bottom: the
    eigenvalues.]{
    \includegraphics[width=6.1cm]{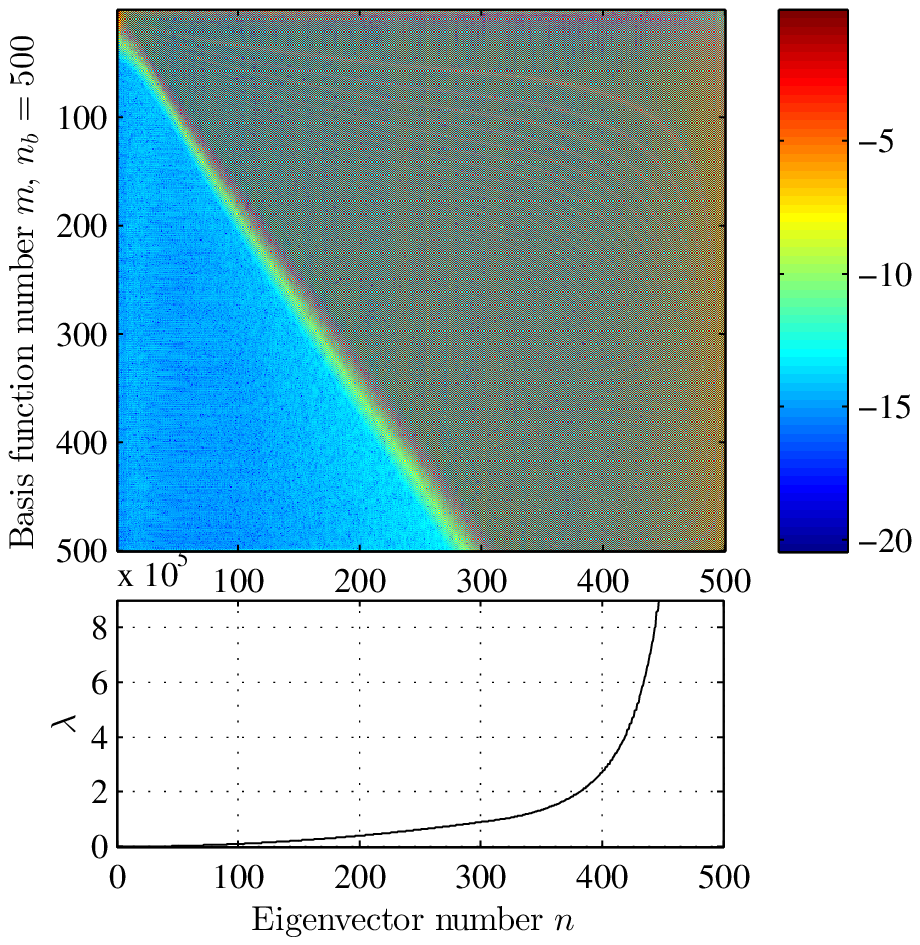}\label{fig:bvpExp2Sub2}}
    \caption{Solution of the Mathieu differential equation. It is important to note that the coefficient $r$ is negative.}
   \label{fig:bvpExp2}
\end{figure}

The second example is a Mathieu differential equation; we have
taken this example from~\cite[Problem 2, with $r=25$]{Pryce1999}.
This equation arises in the vibration of elliptical membranes. We
have chosen this problem because it is known to produce a pair of
closely located eigenvalues; this should enable the test of the
resolution of the eigenvalues computed using the proposed method.
Secondly, the solution to the equation has no known analytical
form, this makes numerical solutions particularly valuable. The
Mathieu differential equation is,
\begin{equation}
    \ddot{y} + 2\,r\,\cos{(2 x)}\, y = \lambda \, y,
    \hspace{5mm}
    \text{with}
    \hspace{5mm}
    y(0) = 0,
    \hspace{2mm}
    \text{and}
    \hspace{2mm}
    y(\pi) = 0
\end{equation}
in the interval $0 \leq x \leq \pi$. A Chebyshev distribution of
$n=1000$ nodes are used covering the complete interval, a support
length $l_s = 13$ and $m=500$ basis functions were used for the
computation. The result of the numerical computation can be seen
in Figures~(\ref{fig:bvpExp2Sub1}) and~(\ref{fig:bvpExp2Sub2}).
The pair of expected eigenvalues are computed as, $\lambda_1 =
-2.131489E+01$ and $\lambda_2 = -2.131486E+01$.

The results of these computations are slightly more difficult to
interpret absolutely since the analytical results are not given.
It can be said that the expected pair of eigenvalues have been
found. Secondly, from the nature of the Rayleigh-Ritz spectrum and
the corresponding computed eigenvalues, see
Figure~(\ref{fig:bvpExp2Sub2}), suggest that approximately the
first $k=350$ eigenvalue are correctly computed.

\subsubsection{Sturm-Liouville Example 3}

This is the truncated hydrogen equation, taken from~\cite[Problem
4]{Pryce1999}. This is an example of s singular Sturm-Liouville
equation with a limit point non-oscillatory (LPN) end point.
Although only one constraint is available the equation is well
conditioned, since $g(x) \rightarrow \infty$ as $x \rightarrow 0$.
Highly accurate results for some eigenvalue are available for this
equation~\cite{Pryce1999}. These values are used to evaluate the
accuracy of the computation method presented here. The
differential equation is,
\begin{equation}
    -\ddot{y} +
    \left(
        \frac{2}{x^2} - \frac{1}{x}
    \right) y = \lambda \, y,
    \hspace{5mm}
    \text{with}
    \hspace{5mm}
    y(0) = LPN,
    \hspace{2mm}
    \text{and}
    \hspace{2mm}
    y(1000) = 0
\end{equation}
in the interval $0 \leq x \leq 1000$.

The computation was performed using $n = 1000$ Chebyshev points on
the range $0 < x < 1000$; further, $m=500$ basis functions were
used, whereby only one constraint is applied, i.e., at $x = 1000$
and a support length of $l_s = 13$. The result of the computation
are presented in Table~(\ref{tab:fourEigs}). The known
eigenvalues, i.e., $\lambda_0$, $\lambda_9$, $\lambda_{17}$ and
$\lambda_{18}$ are comparable up the the $10^{\text{th}}$,
$8^{\text{th}}$, $6^{\text{th}}$ and $5^{\text{th}}$ significant
digits respectively. This indicates a high degree of accuracy,
particularly considering that this is a general framework for
differential equations and not a dedicated Sturm-Liouvalle solver.
In~\cite{Amodio2011} difficulties with oscillations at the right
end point were observed for eigenvalues $\lambda_k$ when $k > 8$,
these difficulties are not observed with the methods proposed
here.
\begin{table}[h]
  \centering
  \begin{tabular}{|l|r|r|r|}
        \hline
         & new method & known value~\cite{Pryce1999}& Rel.
        Error\\
        \hline
        \hline
        $\lambda_0 = $ & $-6.2499999978E-02$ & $-6.2500000000E-02$
        & $3.4874503285E-10$\\
        $\lambda_9 = $ & $-2.0661156136E-03$ & $-2.0661157025E-03$
        & $4.3009091823E-08$\\
        $\lambda_{17} = $ & $-2.5757218232E-04$ & $-2.5757359232E-04$ &
        $5.4741446402E-06$\\
        $\lambda_{18} = $ & $2.8740937561E-05$ & $2.8739013100E-05$ &
        $-6.6963370220E-05$\\
        \hline
  \end{tabular}
  \caption{Table of eigenvalues for the truncated hydrogen equation.
  Comparing the numerical results obtained using the new method with the known results~\cite{Pryce1999}}\label{tab:fourEigs}
\end{table}
%
\subsection{Boundary Value Problem Example 1}
\begin{figure}[h]
    \centering
    \includegraphics[width=6cm,clip=true]{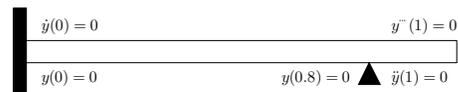}
    \caption{Cantilever with an additional support representing an inner-value problem.}
   \label{fig:BVP1}
\end{figure}

\begin{figure}[h]
    \centering
    \subfigure[The first three admissible functions for the cantilever with and additional simple support.]{
    \includegraphics[width=6.cm]{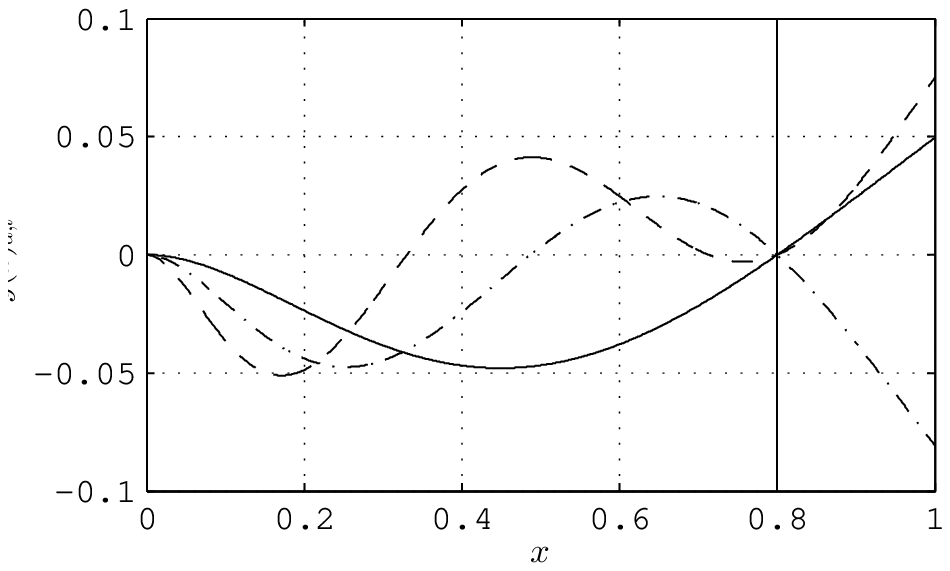}
     \label{fig:BVP1Sub1}}
     \hspace{3mm}
    \subfigure[The first three eigenfunctions.]{
        \includegraphics[width=6cm]{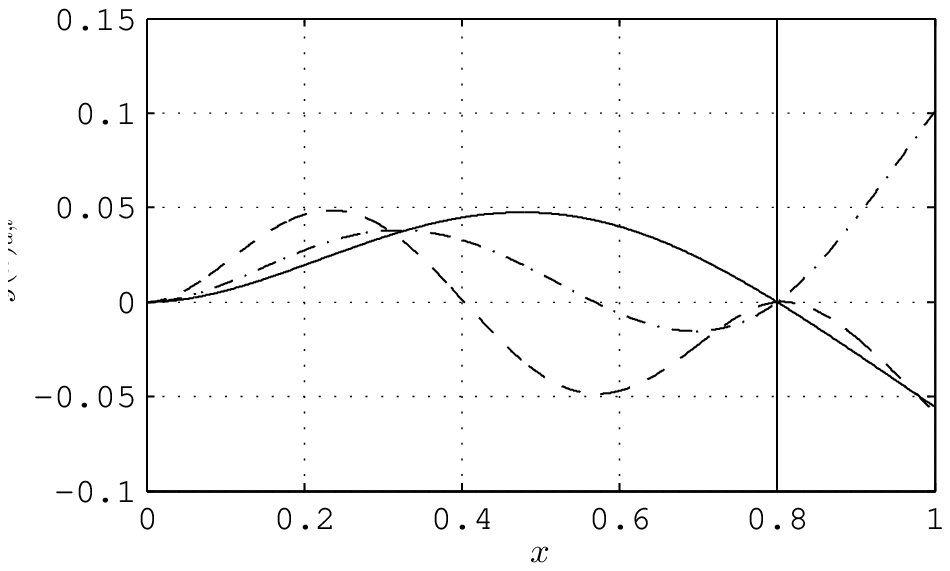}
    \label{fig:BVP1Sub2}}
    \caption{Admissible functions and eigenfunctions for the cantilever with
    additional simple support shown in Figure~(\ref{fig:BVP1}).}
   \label{fig:BVP1Fns}
\end{figure}

The last test is a classical Engineering boundary value problem. A
cantilever with additional simple support forcing the constraint
$y(0.8) = 0$ is shown in Figure~(\ref{fig:BVP1}); this is an
example of an inner constraint. Furthermore, the system is over
constrained, since there are $5$ constraints placed on a
$4^{\text{th}}$ order differential equation. It demonstrates the
ability of the proposed framework to solve problems with
arbitrarily placed constraints and to solve over-constrained
systems. The first two admissible and eigenfunctions are shown in
Figures~(\ref{fig:BVP1Sub1}) and~(\ref{fig:BVP1Sub2})
respectively.
This computation was performed with $n = 1001$ points evenly
distributed along the interval $0 \leq x \leq 1$, a total of $r =
500$ basis functions were used and a support length $l_s = 13$ was
used for the computation of the local derivatives. The method
presented for this problem is a discrete equivalent of a
Rayleigh-Ritz method. The Rayleigh-Ritz coefficients for the first
four eigenfunctions with respect to the first $10$ constrained
basis functions are shown in Table~(\ref{tab:RRCfs}).
\begin{table}[h]
  \centering
  \begin{tabular}{|l|r|r|r|r|}
  \hline
        & $\Phi(x)_1$ & $\Phi(x)_2$ & $\Phi(x)_3$ & $\Phi(x)_4$\\
        \hline
        \hline
  $c_0 $ & -0.99640& -0.06818& -0.04144& -0.00376\\
  $c_1 $ &   0.08434& -0.93235& -0.33425& -0.07821\\
  $c_2 $ &   0.00763& 0.34853& -0.85504& -0.36674\\
  $c_3 $ &   -0.00317& -0.06619& 0.39012& -0.82201\\
  $c_4 $ &   -0.00041& -0.01070& 0.04490& -0.41659\\
  $c_5 $ &   0.00010& 0.00334& -0.02068& 0.04126\\
  $c_6 $ &   -0.00024& -0.00767& 0.02352& -0.08609\\
  $c_7 $ &   0.00016& 0.00522& -0.01475& 0.02897\\
  $c_8 $ &  -0.00005& -0.00172& 0.00487& -0.00615\\
  $c_9 $ &   0.00002& 0.00053& -0.00157& 0.00340\\
  \hline
  \end{tabular}
  \caption{The discrete Raleigh-Ritz coefficients for the first
four eigenfunctions with respect to the first $10$ constrained
basis functions $\M{B}_c$ for the cantilever with additional
simple support (see Figure~(\ref{fig:BVP1})).}\label{tab:RRCfs}
\end{table}
\section{Conclusions}
The successful solution od a series of initial-, boundary- and
inner-value problems, with excellent results, demonstrates the
general applicability of the proposed matrix framework to the
solution of ODEs. The generic formulation pof the solution method
as a least squares problem is very powerful, since it enables the
application of the methods to many classes of problems including
inverse problems.

In the case of initial value problems the least squares solution
ensures that the solution has no implicit direction of solution,
i.e., errors do not accumulate as the computation proceeds. Also
the application of the framework to a selection of Sturm-Liouville
problems has delivered results comparable with those delivered by
dedicated Sturm-Liouville solvers.

The key issues in this paper which led to this success are:
\begin{enumerate}
    \item A Lanczos process with complete
    reorthogonalization is used to synthesize the polynomial basis
    functions. This ensures highly accurate polynomial basis for
    the computation.
    \item A correct definition of the local differentiating matrix
    with consistent degree of
    approximation over the complete support. This ensures the
    possibility of correctly estimating differentials at the
    boundary:  essential for boundary and initial value
    problems.
    \item The formulation of a method of generating orthonormal
    homogeneous admissible functions from constraints. The matrix
    containing these basis functions is ortho-normal, yielding
    optimal behavior in terms of error propagation.
    This enables the implementation of a discrete equivalent of
    the Rayleigh-Ritz method.
    \item The formulation of the solution of the ODE as a least
    squares approximation; this ensure that there is no accumulation of errors.
\end{enumerate}
\bibliographystyle{siam}
\bibliography{BIB/diffMatrices,BIB/diffEq,BIB/CPS}

\end{document}